\documentclass[12pt]{article}
\usepackage{amstex, amsthm}

\def\Gal{\operatorname{Gal}}
\def\FF{\mathbb{F}}
\def\NN{\mathbb{N}}
\def\QQ{\mathbb{Q}}

\def\ZZ{\mathbb{Z}}

\def\Fp{\FF_p}
\def\Fpbar{\overline{\Fp}}
\def\map#1#2#3{{#1\!: #2 \to #3}}

\newtheorem{thm}{Theorem}
\newtheorem{lem}[thm]{Lemma}
\newtheorem{cor}[thm]{Corollary}
\theoremstyle{definition}
\newtheorem{defn}{Definition}

\begin{document}

\title{The algebraic closure of the power
series field in positive characteristic}
\author{Kiran S. Kedlaya \\ Massachusetts Institute of Technology}
\date{December 16, 1999}

\maketitle

\begin{abstract}
For $K$ a field,
let $K((t))$ denote the quotient field of the power series ring over
$K$. The ``Newton-Puiseux theorem'' states that if $K$ has characteristic 0,
the algebraic closure of $K((t))$
is the union of the fields $K((t^{1/n}))$ over $n \in \NN$.
We answer a question of Abhyankar by constructing an algebraic closure
of $K((t))$ for any field $K$ of positive characteristic explicitly in
terms of certain generalized power series.
\end{abstract}

\section{Introduction}

For $K$ a field, let $K((t))$ denote the field of formal power series over
$K$ (that is, expressions of the form $\sum_{i=m}^\infty x_i
t^i$ for some $m \in \ZZ$ and $x_i \in K$, with the usual
arithmetic operations). A classical theorem
\cite[Proposition~II.8]{ser}
attributed to Puiseux, but essentially known to Newton,
states that if $K$ is an algebraically closed field of
characteristic zero, then the algebraic
closure of $K((t))$ is isomorphic to 
\[
\bigcup_{i=1}^\infty K((t^{1/i})).
\]
Hereafter,
we will take $K$ to be an algebraically closed field of characteristic
$p>0$.
In this case,
Chevalley \cite{che} noted that
the Artin-Schreier polynomial $x^p - x - t^{-1}$ has no root in the
Newton-Puiseux
field. In fact, the Newton-Puiseux field is precisely the
perfect closure of the maximal tamely ramified extension of
$K((t))$.

Abhyankar \cite{abh} pointed out that under a suitable generalization
of the notion of power series, Chevalley's polynomial should acquire
the root
\[
x = t^{-1/p} + t^{-1/p^2} + t^{-1/p^3} + \cdots.
\]
The generalization we use here was introduced by Hahn \cite{hahn}, and we
will only give a brief introduction here; for detailed treatments,
see \cite{poo} or \cite{rib}.

A \emph{generalized power series}
(or simply ``series'') is an expression of the form
$\sum_{i \in \QQ} x_i t^i$ with $x_i \in K$,
where the set of $i$ such that $x_i \neq 0$ (the
\textit{support} of the series) is a \emph{well-ordered} subset of $\QQ$,
that is, one every subset of which has a least element.
We add and multiply generalized power series in the natural way:
\begin{align*}
\sum_{i} x_i t^i + \sum_{j} y_j t^j
&= \sum_{k} (x_k + y_k) t^k \\
\sum_{i} x_i t^i \cdot \sum_{j \in T} y_j t^j
&= \sum_{k} \left( \sum_{i+j=k}
 x_i y_j \right) t^k.
\end{align*}
Note that multiplication makes sense because for any $k$, there are only
finitely many pairs $i,j$ with $i+j = k$ and $x_iy_j \neq 0$. Also, 
both the sum and product have well-ordered supports, so the generalized
power series form a ring under these operations.

The ring
of generalized power series is quite large, so one
might reasonably expect it to contain an algebraic closure of
$K((t))$. A stronger assertion was proved independently
by Huang \cite{hua}, Rayner \cite{ray} and \c{S}tef\u anescu \cite{ste}. (Huang's PhD thesis,
written under Abhyankar, does not appear to have been published.)
\begin{thm}[Huang, Rayner, \c{S}tef\u anescu]
Let $L$ be the set of generalized power series of the form
$f = \sum_{i \in S} x_i t^i$ $(x_i \in K)$, where the set $S$ (which 
depends on $f$) has the
following properties:
\begin{enumerate}
\item
Every nonempty subset of $S$ has a least
element (i.e.\ $S$ is well-ordered).
\item
There exists a natural number $m$ such that every element of $mS$ has
denominator a power of $p$.
\end{enumerate}
Then $L$ is an algebraically closed field.
\end{thm}

The purpose of this paper is to refine this result, by
determining precisely which series of the form described by Theorem~1 are algebraic over $K((t))$. Such series will satisfy two additional restrictions: a further condition on the support of the series,
and (unlike in characteristic zero) a condition on the coefficients themselves.
A prototype of the latter condition is the following result, also due independently to Huang and to Stef\u anescu~\cite{ste2}.
\begin{thm}[Huang, Stef\u anescu]
The series $\sum_{i=1}^\infty x_i t^{-1/p^i}$ $ (x_i \in \Fpbar)$
is algebraic over $\Fpbar((t))$ if and only if the sequence $\{x_i\}$ is
eventually periodic.
\end{thm}
Our results imply Theorem~2 as well as
other results of Benhissi~\cite{ben}, Huang, and Vaidya~\cite{vai}. (They do
not directly imply Theorem~1, but our approach can be easily adapted to give a
short proof of that theorem.)

It should be noted that an analogous description of the algebraic closure of
a mixed-characteristic complete discrete valuation ring can be given;
see \cite{me} for details.

\section{Lemmas}

We begin with two preparatory lemmas. The first lemma is a
routine exercise in Galois theory.

\begin{lem} \label{gal}
Every finite normal extension of $K((t))$ is contained in a
tower of Artin-Schreier extensions over $K((t^{1/n}))$ for some $n
\in \NN$.
\end{lem}
\begin{proof}
Let $L$ be a finite normal extensions of $K((t))$ of inseparable degree
$q$; then $L$ is Galois over $K((t^{1/q}))$. We now appeal to
results from \cite[Chap.~IV]{ser} on extensions of complete fields:
\begin{enumerate}
\item
A finite Galois extension of complete fields inducing the trivial
extension on residue fields is totally ramified.
\item
The wild inertia group of such an extension is a $p$-group.
\item
The quotient of the inertia group by the wild inertia group is cyclic
of degree prime to $p$.
\end{enumerate}
Let $M$ be the maximal subextension of
$L$ tamely ramified over
$K((t^{1/q}))$ and let $m$ be the degree of $M$ over $K((t^{1/q}))$.
By Kummer theory (since $K$ contains an $m$-th root of unity),
$M = K((t^{1/q}))(x^{1/m})$ for some $x \in K((t^{1/q}))$, but this
directly implies $M = K((t^{1/qm}))$.

Now $L$ is a $p$-power extension of $M$, so to complete the proof of the
lemma, we need only show that $L$ can be expressed as a tower of
Artin-Schreier extensions over $M$. Since every nontrivial $p$-group has a
nontrivial center, we can find a normal series
\[
\Gal(L/M) = G_0 \rhd G_1 \rhd \cdots \rhd G_n = \{1\}
\]
in which $[G_{i-1}:G_{i}]=p$ for $i = 1, \dots, n$.
The corresponding subfields form a tower of degree $p$ Galois
extensions from $M$ to $L$. Every Galois extension of degree $p$ in
characteristic $p$ is an Artin-Schreier extension (by the
additive version of Hilbert's Satz~90),
completing the proof.
\end{proof}

The second lemma characterizes sequences satisfying the
``linearized recurrence relation'' (hereafter abbreviated LRR)
\begin{equation} \label{twist}
d_0 c_n + d_1 c_{n+1}^p + \cdots + d_k c_{n+k}^{p^k} = 0
\end{equation}
for $n \geq 0$. Of course we may assume $d_k \neq 0$; if we are willing
to neglect the first few terms of the series, we may also assume $d_0 \neq 0$.
\begin{lem}
Let $k$ be a positive integer and let $d_0, \dots, d_k$ be elements
of $K$ with $d_0, d_k \neq 0$.
\begin{enumerate}
\item
The roots of the polynomial $P(x) = d_0 x + d_1 x^p + \cdots + d_k
x^{p^k}$ form a vector space of dimension $k$ over $\Fp$.
\item
Let $z_1, \dots, z_k$ be a basis of the aforementioned vector space.
Then a sequence $\{c_n\}$ satisfies (\ref{twist}) if and only if it has the form
\begin{equation} \label{soleq}
c_n = z_1 \lambda_1^{1/p^n} + \cdots + z_k \lambda_k^{1/p^n}
\end{equation}
for some $\lambda_1, \dots, \lambda_k \in K$.
\end{enumerate}
\end{lem}
\begin{proof}
\ 
\begin{enumerate}
\item
The roots of $P(x)$ form an $\Fp$-vector space because $P(x+y) = P(x)
+ P(y)$ (that is, $P$ is linearized), and the dimension is $k$ because
$P'(x) = d_0$ has no zeroes, so $P(x)$ has distinct roots.
\item
The set of sequences satisfying (\ref{twist}) forms a $K$-vector
space with scalar multiplication given by the formula
\[
\lambda \cdot (c_0, c_1, c_2, \dots) = (c_0 \lambda, c_1 \lambda^{1/p}, c_2 \lambda^{1/p^2}, \dots).
\]
(but not with the usual scalar
multiplication, which will cause some difficulties later).
The dimension of this space is clearly $k$, since $c_0, \dots,
c_{k-1}$ determine the entire sequence.
On the other hand, the sequences satisfying (\ref{soleq}) form a
$k$-dimensional subspace, since the Moore determinant
\[
\det \begin{pmatrix} z_1 & \cdots & z_k \\
z_1^p & \cdots & z_k^p \\
\vdots & & \vdots \\
z_1^{p^{k-1}} & \cdots & z_k^{p^{k-1}}
\end{pmatrix}
\]
is nonzero whenever $z_1, \dots, z_k$ are linearly independent over
$\Fp$.
Thus all solutions of (\ref{twist}) are given by (\ref{soleq}).
\end{enumerate}
\end{proof}

\begin{cor} \label{cor}
If the sequences
$\{c_n\}$ and $\{c'_n\}$ satisfy LRRs with coefficients
$d_0, \dots, d_k$ and $d'_0, \dots, d'_\ell$, respectively, then the sequences
$\{c_n + c'_n\}$ and $\{c_n c'_n\}$ satisfy LRRs
with coefficients depending only on the $d_i$ and $d'_i$.
\end{cor}
\begin{proof}
Let $z_1, \dots, z_k$ and $y_1, \dots, y_\ell$ be $\Fp$-bases for the
roots of the polynomials $d_0 x + \cdots + d_k x^{p^k}$ and $d'_0 x + \cdots +
d'_\ell x^{p^\ell}$, respectively. Then for suitable $\lambda_i$ and $\mu_j$,
\begin{align*}
c_n &= \sum_i z_i \lambda_i^{1/p^n} \\
c'_n &= \sum_j y_j \mu_j^{1/p^n} \\
c_n + c'_n &= \sum_i z_i \lambda_i^{1/p^n} + \sum_j y_j \mu_j^{1/p^n} \\
c_n c'_n &= \sum_{i,j} (z_i y_j) (\lambda_i \mu_j)^{1/p^n}.
\end{align*}
In other words, $\{c_n + c'_n\}$ and $\{c_n c'_n\}$ satisfy LRRs 
whose coefficients are those of the polynomials whose roots comprise the
$\Fp$-vector space spanned by $z_i + y_j$ and $z_i y_j$, respectively. In
particular, these coefficients depend only on the $d_i$ and $d'_j$,
and not on the particular sequences $\{c_n\}$ and $\{c'_n\}$.
\end{proof}

\section{The main result: algebraic series over $K((t))$}

We now construct the sets on which algebraic series are
supported. For $a \in \NN$ and $b,c \geq 0$,
define the set
\[
S_{a, b, c} = \left\{ \frac{1}{a}(n - b_1 p^{-1} - b_2 p^{-2} - \cdots):
n \geq -b, b_i \in \{0, \dots, p-1\}, \sum b_i \leq c \right\}.
\]
Since $S_{a,b,c}$ visibly satisfies the conditions of Theorem~1, any
series supported on $S_{a,b,c}$ belongs to the field $L$ of that theorem.
\begin{thm} \label{thm:mini}
The ring of series supported on $S_{a,b,c}$ for some $a,b,c$ contains
an algebraic closure of $K((t))$.
\end{thm}
We do not give an independent proof of this result, as it
follows immediately from Theorem~\ref{thm:main}, which we prove
directly. Beware that
the ring in Theorem~\ref{thm:mini} is not a field! For example,
the inverse of $x = \sum_{i=1}^\infty x_i t^{-1/p^i}$ is not supported on
any $S_{a,b,c}$ unless $x$ is algebraic over $K((t))$.

As noted earlier, to isolate the algebraic closure of $K((t))$ inside
the ring of generalized power series, it does not suffice to constrain
the support of the series; we must also impose a ``periodicity''
condition on the coefficients. Such a condition should resemble the
criterion of Theorem~2, but with two key differences: it should apply
to an arbitrary field $K$, and it must constrain the coefficients of a
series supported on $S_{a,b,c}$, which cannot be naturally organized
into a single sequence.

The first difference is addressed by Vaidya's generalization of
Theorem~2 \cite[Lemma~4.1.1]{vai}: for $K$ arbitrary, the series
$\sum_{i=1}^\infty x_i t^{-1/p^i}$ is algebraic over $K((t))$ if and
only if the sequence $\{x_i\}$ satisfies an LRR. To address the second
difference, we must impose Vaidya's criterion on many different
sequences of coefficients in a uniform way, so that the criterion
actually forces the series to be algebraic. The following
definition fulfills this demand.

\begin{defn} \label{def1}
Let $T_c = S_{1,0,c} \cap (-1, 0)$. A function $\map{f}{T_c}{K}$ is
\emph{twist-recurrent of order $k$,} for some positive integer $k$,
if there exist $d_0, \dots, d_k \in K$ such that the LRR (\ref{twist}) holds
for any sequence $\{c_n\}$ of the form
\begin{equation} \label{seq}
c_n = f(- b_1 p^{-1} - \cdots - b_{j-1} p^{-j+1} -
p^{-n}(b_j p^{-j} + \cdots) ) \quad (n \geq 0)
\end{equation}
for $j \in \NN$ and $b_1, b_2, \dots \in \{0, \dots, p-1\}$ with $\sum b_i
\leq c$.
\end{defn}
An example may clarify the definition: one such sequence is
\[
f(-1.2021), f(-1.20201), f(-1.202001), f(-1.2020001), \dots,
\]
where the arguments of $f$ are written in base $p$.
If the value of $k$ is not relevant, we simply say $f$ is twist-recurrent.

\begin{defn} \label{def2}
A series $x = \sum x_i t^i$ is \emph{twist-recurrent}
if the following conditions hold:
\begin{enumerate}
\item
There exist $a,b,c \in \NN$ such that $x$ is supported on $S_{a,b,c}$.
\item
For some (any) $a,b,c$ for which $x$ is supported on $S_{a,b,c}$,
and for each integer $m \geq -b$,
the function $\map{f_m}{T_c}{K}$ given by $f_m(z)
= x_{(m+z)/a}$ is twist-recurrent of order $k$ for some $k$.
\item
The functions $f_m$ span a finite-dimensional
vector space over $K$.
\end{enumerate}
\end{defn}
Note that condition 3 implies that the choice of $k$ in the second
condition can be made independently of $m$. It does not imply, however, that
the coefficients $d_0, \dots, d_k$ can be chosen independently of $m$ (except
in the case $K = \Fpbar$, as we shall see in the next section).

\begin{lem} \label{reduc}
Every twist-recurrent series supported on $S_{1,b,c}$ can be written
as a finite $K((t))$-linear combination of twist-recurrent series
supported on $T_c$, and every such linear
combination is twist-recurrent.
\end{lem}
\begin{proof}
If $x = \sum_i x_i t^i$ is a twist-recurrent series supported on
$S_{1,b,c}$, define the functions $\map{f_m}{T_c}{K}$ for $m \geq -b$ by
the formula $f_m(z) = x_{m+z}$, as in Definition~\ref{def2}.
By condition~3 of the definition, the $f_m$ span a finite-dimensional
vector space of functions from $T_c$ to $K$; let $g_1, \dots, g_r$ be a
basis for that space,
and write $f_m = k_{m1}
g_1 + \cdots + k_{mr}g_r$. Now writing
\begin{align*}
x &= \sum_{m \geq -b} \sum_{i \in T_c} f_m(i) t^{m+i} \\
&= \sum_{m \geq -b} \sum_{i \in T_c} t^{m+i} \left( \sum_{j=1}^r k_{mj} g_j(i) \right) \\
&= \sum_{j=1}^r \left( \sum_{m \geq -b} k_{mj} t^m \right) \left( \sum_{i \in T_c} g_j(i) t^i \right)
\end{align*}
expresses $x$ as a finite $K((t))$-linear combination of twist-recurrent
series supported on $T_c$.
Conversely, to show that a linear
combination of twist-recurrent series on $T_c$ is twist-recurrent, it
suffices to observe that by Corollary~\ref{cor},
the sum of two twist-recurrent functions is
twist-recurrent (thus verifying condition~2 for the sum, the
other two being evident).
\end{proof}

\begin{thm} \label{thm:main}
The twist-recurrent series form an algebraic closure of $K((t))$.
\end{thm}

\begin{proof}
We verify the following three assertions.
\begin{enumerate}
\item
Every twist-recurrent series is algebraic over $K((t))$.
\item
The twist-recurrent series are closed under addition and scalar multiplication.
\item
If $y$ is twist-recurrent and $x^p - x = y$, then $x$ is twist-recurrent.
\end{enumerate}
From these, it will follow that the twist-recurrent series form a ring
algebraic over $K((t))$ (which is automatically then a field) closed
unter Artin-Schreier extensions; by Lemma~\ref{gal}, this field is
algebraically closed.

Before proceeding, we note that for each assertion, it suffices to work with
series supported on $S_{a,b,c}$ with $a=1$.

\begin{enumerate}
\item
We proceed by induction on $c$ (with vacuous base case $c=0$);
by Lemma~\ref{reduc}, we need only
consider a series $x = \sum x_i t^i$ supported on $T_c$. Choose $d_0, \dots, d_k$ as in
Definition~\ref{def2}, and let
\[
y = d_0 x^{1/p^k} + d_1 x^{1/p^{k-1}} + \cdots + d_k x.
\]
Clearly $y$ is supported on $T_c$; in fact, we claim it is
supported on $S_{p^{k},0,c-1}$. To be precise, if $j = - \sum
b_i p^{-i}$ belongs to $T_c$ but not to $S_{p^{kc}, 0,
c-1}$, then $\sum b_i = c$, and $b_i = 0$ for $i \leq k$. In particular
$p^k j$ lies in $T_c$, and so
\[
y_j = d_0 x_{p^k j}^{1/p^k} + d_1 x_{p^{k-1}j}^{1/p^{k-1}} + \cdots + d_k x_j = 0
\]
because $x$ is twist-recurrent.

We conclude that $y$ is twist-recurrent (we have just verified condition~1,
condition~2 follows from Corollary~\ref{cor}, and condition~3 is evident).
By the induction hypothesis, $y$ is algebraic over $K((t))$, as
then are $y^{p^k}$ and thus $x$.

\item
Closure under addition follows immediately from Lemma~\ref{reduc};
as for multiplication,
it suffices to show that $xy$ is twist-recurrent whenever
$x = \sum x_i t^i$ and $y = \sum y_i t^i$ are twist-recurrent on $T_c$.
We will prove this by showing that any sequence of the form
\[
c_n = (xy)_{-b_0 - b_1 p^{-1} - \cdots - b_{j-1}p^{-(j-1)}
- p^{-n} (b_{j} p^{-j} +
\cdots)}
\]
becomes, after some initial terms, the sum of a fixed number of pairwise
products of similar sequences derived from $x$ and $y$. Those sequences satisfy
fixed LRRs, so $\{c_n\}$ will as well by Corollary~\ref{cor}.

To verify this claim, recall that $(xy)_k$ is the sum of $x_i y_j$ over all
$i,j \in T_c$ with $i+j=k$. Writing the sum $(-i) + (-j)$ in base $p$,
we notice that for $n$ sufficiently large, there can be no carries across the
``gap'' between $p^{-(j-1)}$ and $p^{-j-n}$. (To be precise,
the sum of the digits of $-k$ equals the sum of the digits of $(-i)$ and $(-j)$
minus $(p-1)$ times the number of carries.) Thus the number of ways to write
$-k$ as $(-i) + (-j)$ is uniformly bounded, and moreover as $k$ runs
through a sequence of indices of the shape in (\ref{seq}), the possible
$i$ and $j$ are constrained to a finite number of similar sequences. This
proves the claim.

\item
Since the map $x \mapsto x^p - x$ is
additive, it suffices to consider the cases when $y$ is supported on
$(-\infty, 0)$ and $(0, \infty)$.

First, suppose $y$ is supported on
$(-\infty, 0) \cap S_{a,b,c}$ for some $a,b,c$; then
\[
x = \sum_i \sum_{n=1}^\infty y_i^{1/p^n} t^{i/p^n} = \sum_i t^i \sum_n
y_{ip^n}^{1/p^n}
\]
is supported on $S_{a,b,b+c}$. We must show that if $-b \leq m \leq 0$,
$b_i \in \{0, \dots, p-1\}$ and $\sum b_i \leq c$, then for any $j$, the sequence
\begin{equation} \label{step3}
c_n = x_{m-b_1 p^{-1} - \cdots - b_{j-1} p^{-(j-1)} - p^{-n}(b_j p^{-j} + \cdots)}
\end{equation}
satisfies a fixed LRR. If $m < 0$ or $j > 0$, then $\{c_n\}$ is the sum
of a bounded number of sequences satisfying fixed LRRs, namely certain
sequences of the $y_i$, so $x$ is twist-recurrent by
Corollary~\ref{cor}. If $m=j=0$, then
\[
c_{n+1}^p - c_n = y_{-b_1 p^{-1} - \cdots - b_{j-1} p^{-(j-1)} - p^{-n}(b_j p^{-j} + \cdots)};
\]
if $\{c_{n+1}^p - c_n\}$
is twist-recurrent with coefficients $d_0, \dots, d_k$, then
$\{ c_n \}$
is twist-recurrent with coefficients $-d_0, d_0 - d_1, \dots, d_k
- d_{k-1}$.

Next, suppose $y$ is supported on $(0, +\infty) \cap S_{a,b,c}$; then
\[
x = -\sum_i \sum_{n=0}^\infty y_i^{p^n} t^{i p^n} =-
\sum_i t^i \sum_n y_{i/p^n}^{p^n}
\]
is also supported on $S_{a,b,c}$.
For $i< p^k$, we have $y_{i/p^n} = 0$ for $n > k+c$, since the first
$c$ fractional digits of $-i/p^n$ in base $p$ will be $p-1$. Thus each
sequence defined by (\ref{step3}) is the sum of a bounded number of
sequences
satisfying fixed LRRs (the exact number and the coefficients
of the LRRs depending on $m$), and so
Corollary~\ref{cor} again implies that $x$ is twist-recurrent.
\end{enumerate}
\end{proof}

\section{Variations}

Having completed the proof of the main theorem, we now formulate some
variations of its statement, all of which follow as easy corollaries. Some of
the modifications can be combined, but to avoid excessive repetition, we
refrain from explicitly stating all possible combinations.

First, we fulfill a promise made in the abstract by describing the algebraic
closure of $L((t))$ where $L$ is an arbitrary perfect
field of characteristic $p$, not necessarily
algebraically closed.
\begin{cor}
Let $L$ be a perfect
(but not necessarily algebraically closed) field of characteristic $p$.
Then the algebraic closure of $L((t))$ consists of all twist-recurrent
series $x = \sum_i x_i t^i$ with $x_i$ in a finite extension of $L$.
\end{cor}
\begin{proof}
The argument given for assertion~1 in the proof of Theorem~\ref{thm:main}
shows that any twist-recurrent series with coefficients in $M$
is algebraic over $M((t))$. To show conversely that any series which is algebraic
over $L((t))$ has coefficients in a finite extension of $L$, let $E$ be a finite
extension of $L((t))$, and $M$ the integral closure of $L$ in $E$. Then a slight
modification of Lemma~3 implies that $E$ can be expressed as a tower of Artin-Schreier
extensions over $M((t^{1/n}))$ for some $n \in \NN$. Now the argument given for
assertion~3 in the proof of Theorem~\ref{thm:main} shows that if $y$ has coefficients in $M$
and $x^p - x = y$, then $y$ has coefficients in $M$ except possibly for its constant coefficient,
which may lie in an Artin-Schreier extension of $M$. We conclude that the coefficients of
any element of $E$ lie in a finite extension of $L$.
\end{proof}

For $L$ not perfect, the situation is more complicated, since if $y$ has coefficients in $M$
and $x^p - x = y$, $x$ may have coefficients which generate inseparable extensions of $M$.
We restrict ourselves to giving a necessary condition for algebraicity in this case.
\begin{cor}
Let $L$ be a field of characteristic $p$, not necessarily perfect.
If $x = \sum_i x_i t^i$ is a generalized power series which is algebraic over $L((t))$,
then the following conditions must hold:
\begin{enumerate}
\item
There exists a finite extension $L'$ of $L$ whose perfect closure contains all of the $x_i$.
\item
For each $i$, let $f_i$ be the smallest nonnegative integer such that $x_i^{p^{f_i}} \in L'$.
Then $f_i - v_p(i)$ is bounded below.
\end{enumerate}
\end{cor}

Next, we explicitly describe the $t$-adic completion of the algebraic closure of
$K((t))$, which occurs more often in practice than its uncompleted counterpart. The
proof is immediate from Theorem~\ref{thm:main}.
\begin{cor}
The completion of the algebraic closure of $K((t))$ consists of all
series $x = \sum_{i \in I} x_i t^i$ such that for every $n \in \NN$,
the series $\sum_{i \in I \cap (-\infty, n)} x_i t^i$ is
twist-recurrent (equivalently, satisfies conditions~1 and~2 of Definition~\ref{def2}).
\end{cor}
From the corollary it follows that truncating an algebraic series (that is, discarding
all coefficients larger than some real number) gives an element in the completion of
the algebraic closure. In fact a stronger statement is true; it appears possible but
complicated to prove this without invoking Theorem~\ref{thm:main}.
\begin{cor}
Let $x = \sum_i x_i t^i$ be a generalized power series which is algebraic over
$K((t))$. Then for any real number $j$, $\sum_{i<j} x_i t^i$ is also algebraic
over $K((t))$.
\end{cor}
\begin{proof}
By the theorem, we may read ``twist-recurrent'' for ``algebraic''. Since
$S_{a,b,c}$ is well-ordered, there is a smallest element of $S_{a,b,c}$
which is greater than $j$, and replacing $j$ by that element reduces us to the case that
$j$ is rational. (In fact it is even the case that
all accumulation points of $S_{a,b,c}$ are rational.)
Now simply note that the definition of twist-recurrence is stable under 
truncation.
\end{proof}

Next, we note that if the field $K$ is endowed with an
absolute value $| \cdot |$,
then for any real number $r$, we can consider the field $K((t))^c$
of power series with positive radius of convergence, which is to say,
if $x = \sum_i x_i t^i$, then $r^i |x_i| \to 0$ as $i \to \infty$
for some $r > 0$.
This definition extends without change to generalized
power series.
\begin{cor}
The algebraic closure of $K((t))^c$ consists of the
twist-recurrent series with positive radius of convergence.
\end{cor}

Finally, we note that specializing Theorem~\ref{thm:main} to the case $K
= \Fpbar$ allows us to simplify its statement using the following definition.
A function $\map{f}{T_c}{\Fpbar}$ is
\emph{periodic of period $M$ after $N$ terms} if
for every sequence $\{c_n\}$ defined as
in (\ref{seq}),
\[
c_{n + M} = c_n \qquad \forall \ n \geq N.
\]
\begin{lem} \label{lem:pd}
A function $\map{f}{T_c}{\Fpbar}$ is twist-recurrent if and only if is
periodic.
\end{lem}
\begin{proof}
The values of a function periodic after $N$ terms, or of a function
twist-recurrent with coefficients $d_0, \dots, d_N$, are determined
by its values on those numbers $i
\in (-1, 0)$ such that no two nonzero digits of $-i$ in base $p$ are
separated by $N$ or more zeroes. In particular, all values of the function
lie in some finite field $\FF_q$.

Now if $f$ is periodic with period $M$ after $N$ terms, then $f$ is
clearly twist-recurrent of order $k = N+M \log_p q$. Conversely, suppose
$f$ is twist-recurrent for some $d_0, \dots, d_k$; then $f$ is
periodic with period at most $q^k$ after $q^k$ terms, since a sequence
satisfying (\ref{twist}) repeats as soon as a subsequence of $k$
consecutive terms repeats.
\end{proof}

Since the set of periodic series for given values of $M$ and $N$ is a
vector space over $\Fpbar$, we can restate condition~3 of Definition~\ref{def2}
as a simple uniformity condition.
\begin{cor}
A series $x = \sum_{i \in S_{a,b,c}} x_i t^i$ is algebraic over
$\Fpbar((t))$ if and only if there exist $M,N \in \NN$ such that for
every integer $m \geq -b$, the function $f(z) = x_{(m+z)/a}$ is
periodic of period $M$ after $N$ terms.
\end{cor}

\section{Desideratum: An Algebraic Proof}

Although the definition of a twist-recurrent function involves infinitely many
conditions, such a function can be specified by a finite
number of coefficients (the exact number depending on $a,b,c$). Thus it
makes sense to ask for an algorithm to compute, for any real number $r$,
enough of the coefficients of a root $\sum_{i \in I} x_i t^i$
of a given polynomial to determine all of the $x_i$ for $i<r$.

In characteristic 0, such an algorithm exists and is well-known:
a standard application of Newton polygons allows one to compute
the lowest-order term of each root, and one can then translate
the roots of the polynomial to eliminate this term in the root of interest,
compute the new lowest-order term, and repeat. This method works because
the supports of the roots are well-behaved: there are only finitely
many nonzero $x_i$ with $i<r$. Since this is not generally true in
characteristic $p$, a different strategy must be adopted. One approach that
works is to determine at the outset a finite set of indices $I$ such that
for any $r$, the values of $x_i$ for $i \in I \cap (-\infty, r)$ determine
the values of $x_j$ for all $j<r$; then for each $i \in I$ in increasing
order, use the
Newton polygon to determine $x_i$, compute $x_j$ for all $j$ less than
the next element of $I$, and translate the polynomial to eliminate all known
low-order terms. It is not hard to see that such a set $I$ exists, though
writing it down explicitly could be somewhat complicated.

While the above argument gives the desired algorithm, it does not
give an algorithmic proof of Theorem~\ref{thm:main},
since we had to assume the theorem to ensure that the set $I$ is finite
and computable. A direct proof (without Galois theory) of the termination
of the algorithm, and hence of Theorem~\ref{thm:main}, would be of
great interest.

\subsection*{Acknowledgments}

The author was supported by an NDSEG Graduate Fellowship.
Thanks to Bjorn Poonen for helpful conversations
(specifically, for suggesting Lemma~\ref{reduc}), and to
the organizers of the 1997 AMS summer conference on
Applications of Curves over Finite Fields, where the author learned
about this problem.

\end{document}